\documentclass[11pt]{amsart}
\newenvironment{PROOF}{\noindent
\textbf{Proof.}}{\hfill{{\sc QED}}\\}
\newenvironment{proofof}[1]{\noindent
               \textbf{#1.}}{\hfill{{\sc QED}}\\}

\newenvironment{Eg}[1]{{\vspace{1 ex}}\noindent {\sc Example.}{#1}
                {\vspace{1 ex}}}

\input xy
 \xyoption{all}

\usepackage{latexsym}
\usepackage{amsmath}
\usepackage{amssymb}
\usepackage{amscd}
\usepackage{epsfig}
\usepackage{multicol}
\usepackage{amsfonts}
\usepackage{amsmath,amsthm}
\usepackage{amssymb}
\usepackage{euscript}
\usepackage{enumerate,calc}
\usepackage{fullpage}

\newtheorem{THM}{Theorem}

\newtheorem{thm}{Theorem}[section]
\newtheorem{theorem}[thm]{Theorem}
\newtheorem{corollary}[thm]{Corollary}
\newtheorem{cor}[thm]{Corollary}
\newtheorem{lemma}[thm]{Lemma}
\newtheorem{prop}[thm]{Proposition}
\newtheorem{definition}[thm]{Definition}
\newtheorem{remark}[thm]{Remark}

\def\to{\rightarrow}
\def\into{\hookrightarrow}

\numberwithin{equation}{section}

\newcommand{\HG}{H^*_{\G}}
\newcommand{\HX}{{\EuScript H}^*\!X}
\newcommand{\BG}{B\G}
\newcommand{\EG}{E\G}

\newcommand{\Zt}{{\Z/2\Z}}

\newcommand{\RPi}{{\R}P^{\infty}}

\newcommand{\im}{im}
\newcommand{\HGX}{\HG(X;\Zt)}
\newcommand{\HGM}{\HG(M;\Zt)}
\newcommand{\HGMG}{\HG(M^{\G};\Zt)}
\newcommand{\G}{G_{\R}}
\renewcommand{\H}{H_{\R}}
\newcommand{\K}{K_{\R}}
\newcommand{\Gs}{\G^*}
\newcommand{\XG}{X^{\G}}

\DeclareMathOperator{\coker}{coker}
\DeclareMathOperator{\Hom}{Hom}

\newcommand{\Z}{\mathbb Z}

\newcommand{\R}{\mathbb R}
\newcommand{\C}{\mathbb C}

\newcommand{\algg}{\ensuremath{\mathfrak{g}}}

\newcommand{\algt}{\ensuremath{\mathfrak{t}}}

\newcommand{\algh}{\ensuremath{\mathfrak{h}}}
\newcommand{\algs}{\ensuremath{\mathfrak{s}}}

\def\iso{\cong}

\def\g{\mathfrak{g}}

\title[Mod 2 cohomology of real loci]{
The mod 2 cohomology of fixed point sets of anti-symplectic involutions
}
\date{ \today }

\author[Biss]{Daniel Biss}
\address[Daniel Biss]{Department of Mathematics\\Massachusetts Institute
of Technology 2-251\\Cambridge MA 02139}
\email{daniel@math.mit.edu}
\author[Guillemin]{Victor W. Guillemin}
\address[Victor W. Guillemin]{Department of Mathematics\\
Massachusetts Institute of Technology\\
Cambridge, MA 02139}
\email{vwg@math.mit.edu}
\author[Holm]{Tara S. Holm}
\address[Tara Holm]{Department of Mathematics \# 3840\\ University of CA, Berkeley 94720-3840}
\email{tsh@math.berkeley.edu}

\begin{document}

\bibliographystyle{plain}

\begin{abstract}
Let $M$ be a compact, connected symplectic manifold with a Hamiltonian
action of a compact $n$-dimensional torus $G=T^n$. Suppose that
$\sigma$ is an anti-symplectic involution compatible with the
$G$-action.  The real locus of $M$ is $X$, the fixed point set of
$\sigma$.  Duistermaat uses Morse theory to give a description of the
ordinary cohomology of $X$ in terms of the cohomology of $M$.  There
is a residual $\G=(\Zt)^n$ action on $X$, and we can use Duistermaat's
result, as well as some general facts about equivariant cohomology, to
prove an equivariant analogue to Duistermaat's theorem. In some cases,
we can also extend theorems of Goresky-Kottwitz-MacPherson and
Goldin-Holm to the real locus.
\end{abstract}

\maketitle

\section{Introduction}\label{se:intro}

Atiyah observed in \cite{atiyah} that if $M$ is a compact
symplectic manifold and $\tau$ a Hamiltonian action of an
$n$-dimensional torus $G$ on $M$, then the cohomology groups of
$M$ can be computed from the cohomology groups of the fixed point
set $M^G$ of $\tau$.  Explicitly, if $\Phi:M\to\algt^*$ is the
moment map for $\tau$, then a generic component $\Phi^\xi$ of
$\Phi$ is a perfect Bott-Morse function.  Using $\Phi^\xi$, we may
compute
\begin{eqnarray}\label{1}
H^{\ast}(M;\R) & = & \sum_{i=1}^NH^{\ast-d_i}(F_i;\R),
\end{eqnarray}
where the $F_i$ are the connected components of $M^G$ and $d_i$ is
the Bott-Morse index of $F_i$.  This result is also true in equivariant
cohomology:
\begin{equation}\label{2}
H^{\ast}_G(M;\R)  =  \sum_{i=1}^NH^{\ast-d_i}_G(F_i;\R)  =
\sum_{i=1}^NH^{\ast-d_i}(F_i\times BG;\R).
\end{equation}
This is a consequence of Atiyah's result and Kirwan's equivariant
formality theorem for Hamiltonian $G$-manifolds, as shown in
\cite{kirwan}.

In \cite{duis}, Duistermaat proved a ``real form'' version of
(\ref{1}).  Let $\sigma:M\to M$ be an anti-symplectic involution with
the property that
\begin{eqnarray}\label{3}
\sigma\circ\tau_g & = & \tau_{g^{-1}}\circ\sigma
\end{eqnarray}
and let $X=M^{\sigma}$ be the fixed point set of $\sigma$.  We call
$X$ the {\em real locus} of $M$.  The motivating example of this setup
is a complex manifold $M$ with a complex conjugation $\sigma$.
Duistermaat proved that
\begin{eqnarray}\label{4}
H^{\ast}(X;\Z/2\Z) & = & \sum_{i\in
I}H^{\ast-\frac{d_i}{2}}(F_i^{\sigma};\Z/2\Z),
\end{eqnarray}
where $I\subseteq \{ 1,\dots,N\}$ is the set for which $\sigma$ 
preserves $F_i$. The $\Z/2\Z$ coefficients are essential
here; the theorem does not hold with real coefficients. (See the
comments at the end of this section.)

The first of the four main theorems of this paper is an equivariant
analogue of (\ref{4}) similar to the equivariant analogue  (\ref{2}) of
Atiyah's result (\ref{1}).  By (\ref{3}), the group
\begin{eqnarray}\label{5}
\G  =  \{ g\in G\ |\ g^2=id\} \iso  (\Zt )^n
\end{eqnarray}
acts on $X$ and we will prove the following theorem in
Section~\ref{se:eqRLocus}.

\begin{THM}\label{th:additive}
Suppose $M$ is a symplectic manifold with a Hamiltonian action
$\tau$ of a torus $T^n=G$ and an anti-symplectic involution $\sigma$.
Let $X=M^{\sigma}$  denote the real locus of $M$.  Then the group $\G$
acts on $X$, and the $\G$-equivariant cohomology of $X$ with $\Zt$
coefficients is
\begin{eqnarray}\label{6}
\HG(X;\Zt) & = & \sum_{i\in I}H_{\G}^{*-\frac{d_i}{2}}
(F_i^{\sigma};\Zt).
\end{eqnarray}
As above, the subset $I\subseteq \{ 1,\dots,N\}$ is the set for
which $\sigma$ preserves $F_i$.
\end{THM}

The idea of the proof will be to derive (\ref{6}) from (\ref{4}) by a
simple trick.

The second of the main theorems concerns the structure
of
$$
\HG (X;\Zt)
$$
as a module over the ring
\begin{eqnarray}\label{8}
\HG=\HG(pt)& = & \Zt[x_1,\dots,x_n].
\end{eqnarray}
Let $\BG$ and $\EG$ denote the classifying space and classifying bundle of
$\G$.  Then by the Borel definition of equivariant cohomology
$$
\HG(X;\Zt)=H^*(X\times_{\G}\EG;\Zt).
$$
The cohomology on the right hand side can be computed by the
spectral sequence associated with the fibration
$$
X\times_{\G}\EG\to\BG,
$$
and we will deduce from this computation the following theorem.
\begin{THM}\label{Tfreemod}
The equivariant cohomology $\HGX$ is a free module over $\HG$ generated
in dimension zero.  Moreover, as an $\HG$ module, $\HGX$ is isomorphic to
\begin{equation}
\HG\otimes_{\Zt}H^*(X;\Zt).
\end{equation}
\end{THM}
The idea of the proof is to show that this spectral sequence
collapses at its $E_2$ term, that is, that $X$ is equivariantly
formal.  In particular, this implies that $\HG(X)$ is a free
module over $\Zt[x_1,\dots,x_n]$.

The isomorphisms (\ref{1}), (\ref{2}), (\ref{4}), and (\ref{6}) are all
isomorphisms in {\em additive} cohomology.  The next two sections of
this paper will be concerned with the ring structure of $\HG(X;\Zt)$.
Henceforth, we will assume that the fixed point set $M^G$ is finite
and that the one-skeleton of $\tau$,
$$
M^{(1)}=\{p\in M\ |\ \dim(G\cdot p)\leq 1\},
$$
is of dimension $2$.  It is not hard to see that these two assumptions
imply that
$$
M^{(1)}=\bigcup_{i=1}^{N}E_i,\mbox{\phantom{boo}} E_i\cong\C P^1.
$$
Moreover, $E_i$ is point-wise fixed by an $(n-1)$-dimensional torus
$H_i$, and the diffeomorphism
$$
E_i\overset{~}{\to}\C P^1
$$
intertwines the action of $G/H_i$ on $E_i$ with the standard $S^1$
action on $\C P^1$.  In particular, $E_i$ contains exactly two
$G$-fixed points.  Thus, one can describe the intersection properties
of the $E_i$ by a graph $\Gamma$ with edges $e_i$ corresponding to the
spheres $E_i$ and vertices $V_{\Gamma}=M^G$.  Two vertices $p$ and $q$
are joined by the edge $e_i$ if $E_i^G=\{ p,q\}$.  Furthermore, each
edge, $e_i$ is labeled by a weight $\tilde{\alpha}_{e_i}$ of $G$, the
weight associated with the intertwining homomorphism
$$
G/H_i\to  S^1.
$$

Let $i:M^G\to M$ be the inclusion of the fixed points into the
manifold, and consider the induced map in equivariant cohomology,
\begin{eqnarray}\label{9}
i^*:H_G^*(M;\C) & \to & H_G^*(M^G;\C).
\end{eqnarray}
By a theorem of Kirwan \cite{kirwan}, this map is injective.
Moreover,
\begin{eqnarray}\label{10}
H_G^*(M^G;\C) & = & \bigoplus_{p\in M^G}H_G^*(\{ p\} ;\C),
\end{eqnarray}
and since
$$
H_G^*(\{ p\} ;\C)\iso S(\algg^*),
$$
one can regard an element of $H_G^*(M^G;\C)$ as a map
\begin{eqnarray}\label{11}
f:M^G\to S(\algg^*).
\end{eqnarray}
Goresky, Kottwitz, and MacPherson \cite{GKM} computed the image of
$i^*$, thus determining not just the additive equivariant
cohomology of $M$, but in fact the ring structure of this
cohomology.
\begin{thm}\cite{GKM}\label{th:GKM}
A map $f:M^G\to S(\algg^*)$ is in the image of $i^*$ if and
only if for each edge $e_i=\{ p,q\}$ of $\Gamma$
\begin{eqnarray}\label{12}
f(p)-f(q)\in\tilde{\alpha}_{e_1}\cdot S(\algg^* ).
\end{eqnarray}
\end{thm}

The third of the four main theorems of this paper will be a $\Zt$
version of this result for the manifold $X$.  We define the
one-skeleton of the real locus to be the set
\begin{eqnarray}\label{real1-skel}
X^{(1)}=\{ x\in X\ |\ \# (\G\cdot x)\leq 2\}.
\end{eqnarray}
Assume in addition to the above that $M^G=X^{\G}$ and the real locus
of the one-skeleton is the same as the one-skeleton of the real
locus. We will call a manifold with these properties a {\em $\!\!\mod
2$ GKM manifold}. The map
\begin{eqnarray}\label{13}
i^*:\HG(X;\Zt) & \to & \HG(X^{\G};\Zt)
\end{eqnarray}
is injective, and by factoring (\ref{13}) through the map
\begin{eqnarray}\label{14}
\HG(X;\Zt) & \to & \HG(E_i^{\sigma};\Zt),
\end{eqnarray}
where $E_i^{\sigma}=E_i\cap X\iso \R P^1$, we obtain the following
theorem.
\begin{THM}\label{th:realGKM}
Suppose $M$ is a $\!\!\mod 2$ GKM manifold. An element
$$
f\in\HG(X^{\G};\Zt)
$$
can be thought of as a map $f:V_{\Gamma}\to \Zt[x_1,\dots,x_n]$, and
such a map $f$ is in the image of $i^*$ if and only if, for each edge
$e=\{ p,q\}$ of $\Gamma$
$$
f_p-f_q\in\alpha_e\cdot\Zt[x_1,\dots,x_n],
$$
where $\alpha_e\in\Zt[x_1,\dots,x_n]$ is the image of the weight
$\tilde{\alpha}_e$.
\end{THM}
This completely determines the ring structure of $\HG(X;\Zt)$.
This theorem is proved independently by Schmid \cite{schmid} using
different techniques.  The main examples of $\!\!\mod 2$ GKM
manifold include the real loci of non-singular projective toric
varieties and real loci of coadjoint orbits, including
Grassmannian and flag varieties.  In the case of toric varieties, this
equivariant cohomology ring was computed already by Davis and
Januszkiewicz \cite{DJ}, but our description is quite different from
theirs.

In \cite{gh}, Goldin and Holm generalize the GKM result to the
case where the one-skeleton has dimension at most $4$.  Assume in
addition to the dimension hypothesis that $M^G=X^{\G}$ and the
real locus of the one-skeleton is the same as the one-skeleton of
the real locus. We will call a manifold with these properties a
{\em $\!\!\mod 2$ GH manifold}. The last of the main theorems is a
$\Zt$ version of the result of Goldin and Holm for the real locus
$X$.  For a subgroup $\H$ of $\G$, we will let $\pi_{\H}^*$ will
denote the change of coefficient map
$$
\pi_{\H}^*:\HG\to H_{\H}^*
$$
associated with the inclusion $\H\into\G$. In Section~\ref{se:realGH},
we will prove the following theorem.
\begin{THM}\label{th:realGH}
Suppose that $M$ is a $\!\!\mod 2$ GH manifold with $G$ fixed points
fixed points $M^G=\{ p_1,\dots ,p_d\}$. Let $f_i\in \HG$
denote the restriction of $f\in H_{\G}^*(X)$ to the fixed point $p_i$.
The image of the injection $i^*:H_{\G}^*(X)\rightarrow
H_{\G}^*(X^{\G})$ is the subalgebra of functions $(f_1,\dots
,f_d)\in\bigoplus_{i=1}^d \HG$ which satisfy
$$
\left\{\begin{array}{ll}
\pi_{\H}^*(f_{i_j})=\pi_{\H}^*(f_{i_k}) & \mbox {if } \{ p_{i_1},\dots
,p_{i_l}\}=Z_{\H}^{\G}  \\
\sum_{j=1}^l\frac{f_{i_j}}{\alpha_1^{i_j}\alpha_2^{i_j}}\in
\HG & \mbox {if }\{ p_{i_1},\dots ,p_{i_l}\}=Z_{\H}^G \mbox{ and }
\dim Z_{\H} =4
\end{array}
\right .
$$
for all subgroups $\H$ of $\G$ of order $|\H|=2^{n-1}$ and all
connected components $Z_{\H}$ of $X^{\H}$, where
$\alpha_1^{i_j}$ and $\alpha_2^{i_j}$ are the (linearly dependent)
weights of the $\G$ action on $T_{p_{i_j}}Z_{\H}$.
\end{THM}

We re-emphasize that Duistermaat's techniques only apply to additive
cohomology.  Since we are able to obtain results concerning the ring
structure of the equivariant cohomology and its relationship to
ordinary cohomology, we also obtain statements about the ring
structure of the ordinary cohomology as well. Indeed, in many cases,
Duistermaat's isomorphism (\ref{4}) turns out to give a ring
isomorphism. (See Corollaries~\ref{co:ring1} and \ref{co:ring2} to
Theorem~\ref{th:realGKM} and Corollaries~\ref{co:ring3} and
\ref{co:ring4} to Theorem~\ref{th:realGH}.)  When describing these ring
isomorphisms, we will make use of the following notation.  The symbol
$$
H^{2*}(M;\Zt)
$$
will denote the subring
$$
\bigoplus_iH^{2i}(M;\Zt)\subseteq H^*(M;\Zt),
$$
endowed with a new grading wherein a class in $H^{2i}(M;\Zt)$ is
given degree $i$ (and similarly for equivariant cohomology).  Then
under suitable hypotheses, the additive isomorphism of Duistermaat
becomes an isomorphism of graded rings.

In Section~\ref{se:strings}, we discuss an application of our main
theorems to string theory. The $\Zt$-equivariant
cohomology ring of $T^n$ with $\Zt$ coefficients
classifies all possible orientifold configurations of Type II string
theories, compactified on $T^n$.  We explain how to compute this
cohomolgy ring.

The last section of the paper contains some applications of these
results to elementary problems in combinatorics.  A typical such
application is the following. Let $\Gamma$ be the permutahedron, the
Cayley graph of the symmetric group $S_n$ with edges generated by
transpositions.  By definition, the vertices of $\Gamma$ are elements
of $S_n$ and two vertices $\sigma$ and $\tau$ are joined by an edge
if $\tau\sigma^{-1}$ is a transposition. Our goal is to attach
to each vertex $\sigma$ a subset $S_{\sigma}$ of $\{ 1,\dots ,n\}$
such that, for all pairs $\sigma$ and $\tau$ of adjacent vertices,
either $S_{\sigma}=S_{\tau}$ or the symmetric difference
$$
(S_{\sigma}-S_{\tau})\cup(S_{\tau}-S_{\sigma})
$$
is $\{ i,j\}$, where $\tau\sigma^{-1}=(ij)$.

Let $\mathcal{F}_{n+1}$ be the real flag variety in $n+1$ dimensions.
We will prove that the set of solutions to this problem can be
identified with the set
$$
\HG(\mathcal{F}_{n+1};\Zt),
$$
where $\G={\Zt\times\cdots\Zt}$ is the $n$-fold product of $\Zt$.
The results in
this section were inspired by a remark of Ethan Bolker, who pointed
out to us that $\Zt$ representation theory is simply Boolean algebra.

We conclude these prefatory remarks with a few comments about the
$\Zt$ coefficients.  We recall Witten's recipe for computing the
homology of compact manifolds by Morse theory.  Let $X$ be a compact
manifold, $f:X\to\R$ a Morse function, and $C_f^i$ the index $i$
critical set of $f$.  Let $\mathcal{C}_i$ be the vector space
$$
\mathcal{C}_i=\bigoplus_{p_i\in C_f^i} p_i\R
$$
with basis $C_f^i$.

Equip $X$ with a Riemannian metric and let $v$ be the gradient vector
field of $f$.  For generic metrics, the stable and unstable manifolds
of $v$ intersect transversally.  In particular, for every critical
point $p\in C_f^i$, there are a finite number of gradient curves
joining $p$ to critical points
\begin{eqnarray}\label{15}
q_1,\dots ,q_m\in C_f^{i-1}.
\end{eqnarray}
Moreover, each of these points can be assigned an intrinsic
orientation $\varepsilon (p,q_j)\in\{\pm 1\}$.  Now define a boundary
operator
$$
\partial :\mathcal{C}_i\to\mathcal{C}_{i-1}
$$
by setting
$$
\partial p=\sum \varepsilon(p,q_j)\cdot q_j.
$$
Witten \cite{witten} was the first to explicitly formulate Morse
theory in this way; he showed that $\partial$ {\em is} a boundary
operator, namely $\partial^2=0$, and that $H_*(X;\R)$ is the homology
of the complex $(\mathcal{C},\partial)$.

In particular, when all the critical points are of even index,
$\partial$ is automatically zero.  Thus, one gets
\begin{eqnarray}\label{16}
\dim(H_i(X;\R)) & = & \left\{\begin{array}{lr}0, &  i\mbox{ odd,}\\
                    \#\{p\in C_f^i\}, & i\mbox{
even.}
                \end{array}\right.
\end{eqnarray}
This fact is key to Atiyah's result (\ref{1}).  He
observes that if $\# M^G<\infty$ and if $f$ is a generic component of
the moment map, then $f$ is a Morse function with critical points all
of even index, so (\ref{1}) is a special case of (\ref{16}).

This recipe for computing homology also works in characteristic two;
however, when $\Zt$ symmetries are present, the gradient curves
joining the points $p\in C_f^i$ to the points on the list (\ref{15})
often occur in pairs.  For instance, for manifolds which satisfy the
GKM hypotheses, these pairs of curves correspond to the edges $e_i$ of
the graph $\Gamma$.  That is, each $E_i^{\sigma}\iso \R P^1$ in
(\ref{14}) contains of a pair of gradient curves joining the two
vertices of $e_i$. Hence, the $\!\!\mod 2$ version of $\partial$ is
identically zero.

One can obtain Duistermaat's result by exploiting this phenomenon.
The goal of this article is to push these ideas forward by
systematically applying these techniques to the equivariant setting.

\noindent {\bf Acknowledgments}.  The authors would like to thank
the anonymous referees for useful comments and suggestions.

\section{The equivariant cohomology of the real locus}\label{se:eqRLocus}

Recall $M^{2d}$ is a symplectic manifold with a Hamiltonian
action $\tau$ of a torus $G=T^n$.  Suppose further that there is an
anti-symplectic involution $\sigma:M\to M$ with the property that
$$
\sigma\circ\tau_g  = \tau_{g^{-1}}\circ\sigma.
$$
Let $X=M^{\sigma}$ be the fixed point set of $\sigma$.  We call $X$
the {\em real locus} of $M$.  Recall that
Duistermaat proved the following equality, computing the ordinary
cohomology of the real locus.
\begin{eqnarray}\label{duis}
H^{\ast}(X;\Z/2\Z)  = \sum_{i\in
I}H^{\ast-\frac{d_i}{2}}(F_i^{\sigma};\Z/2\Z),
\end{eqnarray}
where $I\subseteq \{ 1,\dots,N\}$ is the set for which $\sigma$ 
preserves $F_i$, and the $d_i$ are the indices of the fixed
point sets $F_i$.

We will prove the equivariant analogue (\ref{6}) to this equality,
computing the additive structure of the equivariant cohomology:
$$
\HG(X;\Zt)  =  \sum_{i\in I}H_{\G}^{*-\frac{d_i}{2}}
(F_i^{\sigma};\Zt).
$$

\begin{proofof}{Proof of Theorem~\ref{th:additive}}
Consider the product action of $T^n$ on
$$
M\times\underbrace{(\C^d\times\cdots\times\C^d)}_n,
$$
in which each $S^1$ factor acts by multiplication on the corresponding
factor of $\C^d$. This is a Hamiltonian action.  If $(\phi_1,\dots
,\phi_n)=\Phi:M\to\R^n$ is the moment map  associated with $\tau$,
then the moment map of this product action is $\Psi=
(\psi_1,\dots\psi_n)$, with
$$
\psi_i(m,z_{1,1},\dots,z_{1,d},\dots,z_{d,d})=\phi_i(m)+
\sum_{j=1}^d|z_{i,j}|^2.
$$
Let $a=(a_1,\dots a_n)\in\R^n$.  If $a_i>\sup(\phi_i)$ for every
$i$, then $\Psi^{-1}(a)$ and $M\times
S^{2d-1}\times\dots\times S^{2d-1}$ are
equivariantly diffeomorphic, so the reduced space
$$
M_{\mbox{{\footnotesize red}}}=M/\!/_aT^n=\psi^{-1}(a)/T^n
$$
is diffeomorphic to $M\times_{T^n}(S^{2d-1}\times\cdots\times
S^{2d-1})$.  Moreover, there is another action of $T^n$ on
$M\times\C^d\times\cdots\times\C^d$, namely $\tau$ coupled with the
trivial action on $(\C^d)^n$. Since this commutes with the product
action, it induces a Hamiltonian action of $T^n$ on
$M_{\mbox{{\footnotesize red}}}$. In addition, one gets from $\sigma$
an involution
$$
(m,z_1,\dots,z_d)\mapsto (\sigma(m),\overline{z_1},\dots,\overline{z_d})
$$
of $M\times\C^d\times\cdots\times\C^d$.  This induces an
anti-symplectic involution $\tilde{\sigma}$ on
$M_{\mbox{{\footnotesize red}}}$.  Thus, one can apply Duistermaat's
theorem to $M_{\mbox{{\footnotesize red}}}$ to get a formula for the
cohomology of the space
$$
M_{\mbox{{\footnotesize red}}}^{\tilde{\sigma}}=
X\times_{\G}
(S^{d-1}\times\cdots\times S^{d-1})
$$
in terms of the cohomology of the spaces
$$
Z_i^d:=
F_i^{\sigma}\times_{\G} (S^{d-1}\times\cdots\times S^{d-1}) =
F_i^{\sigma}\times (\R P^{d-1}\times\cdots\times \R P^{d-1}).
$$
Now $F_i^{\sigma}\times B\G$ is obtained from $Z_i^d$ by
attaching cells of dimension $d$ and higher. So, for fixed $k$, the
sequence $H^k(Z_i^d;\Zt)$ stabilizes
as $d$ grows large, and moreover is equal to the equivariant
cohomology of $X$. Thus one obtains from (\ref{duis}) the
following real analogue:
$$
\HG(X;\Zt)=\sum H_{\G}^{*-\frac{d_i}{2}}(F_i^{\sigma}; \Zt),
$$
where $\G=\Zt\times\cdots\times\Zt$.
\end{proofof}

\section{A spectral sequence}\label{se:specseq}

The goal of this section is to determine the structure of $\HGX$ as a
module over the ring $\HG=\Zt[x_1,\dots,x_n]$.  We do this by calculating the
$E_2$-term of the Leray-Serre spectral sequence converging to $\HGX$
and deducing,
by dimensional considerations obtained in the previous section, that
the spectral sequence must collapse.  This then gives us the desired
statement about the $E_{\infty}$-term and hence $\HGX$.

Recall that by definition, we have
$$\HGX = H^*(X\times_{\G}\EG;\Zt)$$
where $\EG$ is the total space of the universal $\G$-bundle.
Denote this fiber product by $E$.  We then have a fibration
$p:E\rightarrow \BG$ with fiber $X$.  Let $\HX$ denote the local
coefficient system on $\BG$ associated to this fibration.  The
$E_2$-term we would like to compute is then
$$E_2 = H^*(\BG;\HX).$$
The computation takes place in two steps.  The first step, which
is the technical heart of the argument, consists of carrying out
the computation in the 1-dimensional case.  The remainder of the
proof consists of a relatively straightforward exercise in
bookkeeping.

\begin{lemma} \label{Lspect}
Let $\G=\Zt$, so that $\BG=K(\Zt,1) = \RPi$.  Then $H^*(\BG;\HX)$ is
generated over $\HG=(\Zt)[x]$ in degree zero by $H^*(X)^{\G}$.
Moreover, the only relation is given by $x\cdot(\alpha+\nu(\alpha)) =
0$ for $\alpha\in H^*(X)$ and $\nu\in\G.$
\end{lemma}

\begin{PROOF}
By definition, the cohomology $H^*(\BG;\HX)$ that we would like to
compute is the group cohomology $H^*(\G;H^*(X))$ with respect to the
natural action of $\G$ on $H^*(X);$ we will henceforth pass back and
forth between these two notations without comment.   Our goal is to
understand the cohomology as a
module over the cohomology ring $H^*(\G;\Zt)$ corresponding to the
trivial action of $\G$ on $\Zt.$

Denote the non-trivial element of $\Zt$ by $\nu.$  Consider the
$\G$-module $H^*(X)\oplus H^*(X)$ with $\G$-action defined by the
equation $\nu(\alpha,\beta) = (\nu(\beta),\nu(\alpha))$ for all
$\alpha,\beta\in H^*(X).$  We then get a short exact sequence of
$\G$-modules
\begin{equation}\label{eq:SES}
0\rightarrow H^*(X)\stackrel{f}{\longrightarrow} H^*(X)\oplus
H^*(X) \stackrel{g}{\longrightarrow} H^*(X)\rightarrow 0
\end{equation}
where $f(\alpha) = (\alpha,\alpha)$ and $g(\alpha,\beta) =
\alpha+\beta$.   Of course, it is completely essential to this
identification of the cokernel of $f$ that we
work over $\Zt$; in order for this sequence to be exact, we
must have $g(\alpha,\beta) = \alpha-\beta$, but this would not be a
map of $\Zt$-modules if we did not also have $g(\alpha,\beta) =
\alpha+\beta$.

We would like to consider the long exact cohomology sequence
associated with (\ref{eq:SES}).  Observe that we have isomorphisms
\begin{eqnarray*}
H^*(X)\oplus H^*(X) &\cong &\Zt[\G]\otimes_{\Zt} H^*(X)\\
&\cong & \Z[\G]\otimes_{\Z} H^*(X) \\
&\cong &\Hom_{\Z}(\Z[\G],H^*(X))
\end{eqnarray*}
of $\G$-modules. The first isomorphism is more or less the
definition of the left-hand side, the second follows from the fact
that $H^*(X)$ is 2-torsion, and the third follows from the
finiteness of $\G$.  Thus, the module $H^*(X)\oplus H^*(X)$ is
co-induced and its higher cohomology vanishes.  Moreover, the map
$\alpha\mapsto (\alpha,\nu(\alpha))$ provides an isomorphism of
$H^*(X)$ with $(H^*(X)\oplus H^*(X))^{\G},$ so we have
$$H^0(\G,H^*(X)\oplus H^*(X)) \cong H^*(X)$$ and, of course,
$H^0(\G,H^*(X)) = H^*(X)^{\G}.$ The long exact sequence in
question therefore takes the following form.
$${\footnotesize{\xymatrix{
{0}\ar[r]  &  {H^*(X)^{\G}}\ar[r] & {H^*(X)} \ar[r] & {H^*(X)^{\G}}
\ar `r[d]`[l] `[llld]`[d][dll]
& \\
  & {H^1(\G;H^*(X))} \ar[r] & {0} \ar[r]
& {H^1(\G;H^*(X))} \ar `r[d]`[l] `[llld]`[d][dll] &  \\
  & {H^2(\G;H^*(X))}\ar[r]  & {0} \ar[r] & \dots &
}}}$$

The first map in this sequence is the natural inclusion, and the
  second sends $\alpha$ to $\alpha+\nu(\alpha).$  Therefore,
  $H^1(\BG;\HX)$ is the sought-after quotient of $H^*(X)^{\G}$ by the
  subgroup of all elements of the form $\alpha + \nu(\alpha)$; as a
  result, so is each
  $H^n(\BG;\HX)$ for $n\geq 1.$ Recall that $H^*(\BG) = (\Zt)[x]$
  where $x$ is a class of  degree 1.  The desired result then follows
  from the fact that the connecting homomorphisms are multiplication by $x.$
\end{PROOF}

The remainder of the game consists in playing the results of
Lemma~\ref{Lspect} and Theorem~\ref{th:additive} off of one another.

\begin{corollary}\label{Cintermediate}
The action of $\G$ on $H^*(X)$ is trivial.  Thus, in the case that
$\G=\Zt,$ Theorem~\ref{Tfreemod} holds.
\end{corollary}

\begin{PROOF}
Lemma~\ref{Lspect} computes the $E_2$-term
of the Leray-Serre spectral sequence converging to $\HGX$ in the case
$\G=\Zt.$  Now, the dimensions of the graded pieces of this $E_2$-term
are maximized precisely when $\G$ acts trivially on $H^*(X).$
Moreover, the results of the previous section tell us that
the graded pieces of $\HGX$ have exactly these maximal dimensions.
Since the $E_\infty$-term of a spectral sequence can only be as large
as its $E_2$-term, this tells us that the action must be trivial and
further that
$$
\HGX=E_{\infty}=E_2.
$$
In this case, $H^*(X)^{\G} = H^*(X),$ and the relation
$x\cdot(\alpha+\nu(\alpha)) = 0$ is automatically satisfied, so
Lemma~\ref{Lspect}
tells us that $H^*(\BG,\HX)$ is a free module over $\Zt[x]$ generated
in degree zero by $H^*(X).$  This completes the proof of
Theorem~\ref{Tfreemod} in the case $\G=\Zt.$

The triviality of the action of $\G$ on $H^*(X)$ for
higher-dimensional $\G$ follows by restricting to arbitrary
one-dimensional subtori.
\end{PROOF}

The fact that $\G$ acts trivially on $H^*(X)$ is not new; it can also
be derived from Duistermaat's original argument.  Indeed,
Duistermaat's isomorphism can be seen to be $\G$-equivariant, and the
result then follows from the connectedness of the torus.  Using this
fact would have somewhat simplified our argument, but we chose to give
the above proof so as to avoid appealing to unpublished modifications
of the literature.

We now have the technical input to handle the general case.

\begin{lemma}\label{Lgeneral}
Let $\G\cong(\Zt)^n$, for any positive integer $n$.  Then
$H^*(\BG;\HX)$ is a free module over $\HG=(\Zt)[x_1,\dots,x_n]$
generated by $H^0(\BG;\HG)\cong H^*(X)$.
\end{lemma}

\begin{PROOF}
This is now completely classical.  By Corollary~\ref{Cintermediate},
the action of $\G$ on $H^*(X)$ is trivial, so that $\HX$ is actually
the constant sheaf $H^*(X)$.  Therefore, we have
$$H^*(\BG;\HX) = H^*(\BG;\Zt)\otimes_{\Zt} H^*(X) =
\Zt[x_1,\dots,x_n]\otimes_{\Zt} H^*(X).$$
This completes the proof.
\end{PROOF}

Our goal is now entirely within reach; we need only combine the
results we have proven so far as in the proof of
Corollary~\ref{Cintermediate} to establish the $\HG$-module isomorphism
\begin{equation}\label{eq:module}
\HGX\iso\HG\otimes_{\Zt}H^*(X;\Zt).
\end{equation}

\begin{proofof}{Proof of Theorem~\ref{Tfreemod}}
Lemma~\ref{Lgeneral} tells us that the $E_2$-term
of the Leray-Serre spectral sequence converging to $\HGX$ takes
precisely the form that $\HGX$ itself is asserted to have.  However,
the results of the previous section tell us that the graded pieces of
$\HGX$ have the same dimension as those of this $E_2$-term, and
hence that the spectral sequence collapses.  Therefore,
$$
\HGX=E_{\infty}=E_2,
$$
so $\HGX$ is a free module over $\HG$ generated
in dimension zero.  Thus, its additive structure is as given by
(\ref{eq:module}).
\end{proofof}

\section{Chang-Skjelbred in $\Zt$}\label{se:cs}

As a result of the collapse of the spectral sequence proved in the
previous section, the map
$$
i^*:H_{\G}^*(X;\Zt)\to H_{\G}^*(X^{\G};\Zt)
$$
is an injection.  In the case of the original manifold $M$, the
Chang-Skjelbred theorem \cite{CS} identifies the image of this map.
We prove a $\Zt$ version of that theorem here.

As usual, let $\G=(\Zt)^n$ be the $n$-dimensional ``real torus,''
$M$ a $\G$-manifold, and $M^{\G}$ the fixed point set of the
action; we let $i:M^{\G}\hookrightarrow M$ denote the
inclusion.

\begin{thm}\label{Tchang}
Suppose that $\HGM$ is a free $\HG$-module.  For a subgroup
$\H<\Gs$, let $i_{\H}:M^{\Gs}\hookrightarrow M^{\H}$ denote the
inclusion.  Then we have
$$i^*\HGM = \bigcap_{\substack{\H<
\Gs\\|\H|=2^{n-1}}}i_{\H}^*\HG(M^{\H};\Zt).$$
\end{thm}

Our proof closely models the argument given in \cite{GS}, with
appropriate modifications.

First of all, recall that
$$
\HG=\H^*((\RPi)^n;\Zt) =
\Zt[x_1,\dots,x_n],
$$
with $\deg(x_i)=1$.  Moreover, we may view each $x_i$
as a linear functional $x_i:\G\rightarrow\Zt$, that is, an element
of $\Gs$.  This allows us to identify $\HG$ with the symmetric algebra
$S(\Gs)$, a fact analogous to the ordinary identification $H^*_G\cong
S({\mathfrak g}^*)$ when $G$ is a torus with Lie algebra ${\mathfrak
g}$.  This allows us to view elements of $\HG$ as
polynomial functions on $\G$.

\begin{lemma}\label{LkillK}
Let $\K<\G$ be a subgroup and $\phi:M\rightarrow \G/\K$ a
$\G$-equivariant map.  If $\rho\in S(\Gs)$ annihilates $\K$, then it
must also kill $\HGM$.
\end{lemma}

\begin{PROOF}
We have a sequence of $\G$-equivariant maps
\begin{equation}\label{eq:sequence}
M\rightarrow \G/\K\rightarrow pt.
\end{equation}
Notice that
\begin{eqnarray*}
\HG(\G/\K) & = & H^*(G/\K\times_G EG;\Zt)\\
           & = & H^*(B\K;\Zt) \\
           & = & S(\K^*).
\end{eqnarray*}
The sequence (\ref{eq:sequence}) of $\G$-spaces therefore gives rise to the
following diagram of algebras.
$$\xymatrix{
{\HGM} & {\HG(\G/\K)} \ar@{->}[l] &  {\HG}\ar@{->}[l] \\
 & S(\K^*) \ar@{=}[u] & S(\Gs)\ar@{=}[u]
}$$
Therefore, the map $S(\Gs)\rightarrow\HGM$ defining the module
structure factors as $$S(\Gs)\rightarrow S(\K^*)\rightarrow\HGM,$$ and
the proof is complete.
\end{PROOF}

This observation furnishes us with the fundamental tool in proving
localization theorems for equivariant cohomology.

\begin{prop}\label{Pcomplement}
Let $X$ be a closed $\G$-invariant submanifold of $M$.  For some
positive integer $L$, there exist subgroups $(\K)_1,\dots,(\K)_L$
of $G$, each of which is an isotropy subgroup of some point $p\in
M\backslash X$, such that for any $\alpha_1,\dots,\alpha_L\in\Gs$
with $\alpha_i|_{(\K)_i}=0$, the product
$\alpha_1\alpha_2^2\cdots\alpha_L^2\in\HG$ kills $\HG(M\backslash
X)$.
\end{prop}

\begin{PROOF}
Let $U$ be a $\G$-invariant tubular neighborhood of $X$; it suffices
to prove the desired result for the module $\HG(M\backslash U;\Zt)$.
Now, given any orbit $X_i$ in $M\backslash U$ with isotropy subgroup
$(\K)_i$, we may find a $\G$-invariant open neighborhood $U_i$ of $X_i$
admitting a $\G$-equivariant map $U_i\rightarrow\G/(\K)_i$.  By
compactness, we may cover $M\backslash U$ by finitely many such sets
$U_1,\dots,U_L$.  We now show by induction that for all $r\leq L$, if
$\alpha_1,\dots,\alpha_r$ are elements of $\Gs$ with
$\alpha_i|_{(\K)_i}=0$, then $\alpha_1,\alpha_2^2,\dots,\alpha_r^2$
annihilates $\HG(U_1\cup\dots\cup U_r;\Zt)$.

The case $r=1$ is simply a restatement of Lemma~\ref{LkillK}.  For the
inductive step, consider the Meyer-Vietoris sequence associated to the
cover $U_1\cup\dots\cup U_r = (U_1\cup\dots\cup U_{r-1})\cup U_r$.
Denoting $U_1\cup\dots\cup U_{r-1}$ by $V$, we find an exact sequence
$$H^k_{\G}(V\cap U_r)\longrightarrow H^{k+1}(V\cup U_r)
\longrightarrow H^{k+1}_{\G}(V)\oplus H^{k+1}_{\G}(U_r).$$
Now, since $V\cap U_r\subset U_r$, we have a $\G$-equivariant map
$V\cap U_r\rightarrow X_r$, and so the left-hand term of the sequence
is killed by $\alpha_r$.  Meanwhile, by induction, the right-hand term
is annihilated by
$\alpha_1\alpha_2^2\cdots\alpha_{r-1}^2\alpha_r$, and so the
product $\alpha_1\alpha_2^2\cdots\alpha_r^2$ kills the middle term.
\end{PROOF}

We will also need a relative version of the same result.

\begin{prop}\label{Prelative}
Under the hypotheses of Proposition~\ref{Pcomplement}, the module $\HG(M,X)$
is annihilated by the element
$\alpha_1^2\alpha_2^4\cdots\alpha_L^4\in\HG$.
\end{prop}

\begin{PROOF}
Of course, the map of pairs $(M,X)\rightarrow (M/X,X)$ is an
equivalence, so it suffices to compute $\HG(M/X,X)$. Once again,
let $U$ be a $\G$-equivariant tubular neighborhood of $X$.
We cover $M/X$ by two open sets $U/X$ and $M\backslash X$; since, the
projection map $U/X\rightarrow X/X$ is an equivalence, so we may
identify $\HG(M,X)$ with the kernel of the map
$\HG(M)\rightarrow \HG(U)$. Now, let us write the
Meyer-Vietoris sequence for this cover.
$$H^{k-1}_{\G}(U\backslash X)\longrightarrow
H^{k}_{\G}(M/X)\longrightarrow
H^{k}_{\G}(U/X)\oplus H^k_{\G}(M\backslash X)$$
But by the above discussion, this gives rise to the following exact
sequence.
$$H^{k-1}_{\G}(U\backslash X)\longrightarrow
H^{k}_{\G}(M,X)\longrightarrow H^k_{\G}(M\backslash X)$$
Now, Proposition~\ref{Pcomplement} applies to both
ends of this sequence, so the middle term is killed by
$(\alpha_1\alpha_2^2\cdots\alpha_L^2)^2=\alpha_1^2\alpha_2^4\cdots\alpha_L^4$.
\end{PROOF}

Proposition~\ref{Prelative} gives us the basic localization
results we will need.

\begin{cor}\label{Cloc}
In the setting of Proposition~\ref{Pcomplement}, the kernel and
cokernel of the map $$i^*:\HGM\rightarrow\HG(X;\Zt)$$ are
annihilated by the same element
$\alpha_1^2\alpha_2^4\dots\alpha_L^4\in\HG$.
\end{cor}

\begin{PROOF}
Simply apply Proposition~\ref{Prelative} to the
exact sequence
$$H^k_{\G}(M,X)\longrightarrow H^k_{\G}(M)\longrightarrow H^k_{\G}(X)
\longrightarrow H^{k+1}_{\G}(M,X).$$
\end{PROOF}

\begin{cor}\label{co:inject}
The kernel of the map $i^*:H^*_{\G}(M)\rightarrow
H^*_{\G}(M^{\G})$ is torsion, and hence trivial when $\HGM$ is a free
module.
\end{cor}

\begin{cor}\label{killcok}
For every subgroup $\K<\G$, there exists a monomial
$p=\alpha_1\dots\alpha_N$ annihilating the cokernel of the map
$i^*:H^*_{\G}(M)\rightarrow H^*_{\G}(M^\K)$ such that no $\alpha_i$
vanishes on $\K$.
\end{cor}

\begin{PROOF}
Let $q\in M\backslash M^\K$, and let $\K'$ be its isotropy
subgroup.  Since obviously $\K'\not\supset \K$, there is an
$\alpha':\G\rightarrow\Zt$ with $\alpha'|_{\K'}=0$ and
$\alpha'|_{\K}\neq 0$.  But by Corollary~\ref{Cloc}, we can find a
monomial which is a product of elements of the form $\alpha'$ and
annihilates $\coker i^*$.
\end{PROOF}

Finally, we are ready to prove our main result.

\begin{proofof}{Proof of Theorem~\ref{Tchang}}
First of all, since the map $i$ factors as
$$M^{\G}\stackrel{i_{\H}}{\hookrightarrow}M^{\H}\hookrightarrow M,$$
we know
that for all $\H$, the inclusion $\im(i^*)\subset\im(i_{\H}^*)$
holds.

For the other direction, recall first that $\HGM$ is free; therefore,
by Corollary~\ref{co:inject}
the map $i^*$ is injective and we may consequently view $\HGM$ as
a submodule of $\HGMG$.  Suppose $\{ e_1,\dots,e_k\}$ is an
$S(\Gs)$-basis for $\HGM$.  By Corollary~\ref{killcok}, there is a
monomial $p=\alpha_1\cdots\alpha_N$ with $pe\in\HGM$ for every
$e\in\HGMG$.  Thus, we may write
$$pe=f_1e_1+\dots f_ke_k$$
for unique $f_i\in S(\Gs)$.  Now, since $S(\Gs)$ is a unique
factorization domain, we may divide both sides of this identity by $p$
and cancel common factors to obtain the formula
\begin{equation}\label{factor}
e=\frac{g_1}{p_1}e_1+\dots+\frac{g_k}{p_k}e_k
\end{equation}
where the $g_i$ are uniquely-determined elements of $S(\Gs)$ and the
$p_i$ are uniquely-determined divisors of $p$ such that $g_i$ and
$p_i$ are relatively prime.

Suppose now that $e$ were actually in $\im(i_{\H}^*)$.  We may
find a subset $\{j_1,\dots,j_R\}\subset\{1,\dots,N\}$ such that no
$\alpha_{j_i}$ kills $\H$, and $q=\alpha_{j_1}\cdots\alpha_{j_R}$
annihilates the cokernel of the map
$$
\HGM\rightarrow\HG(M^{\H};\Zt).
$$
Therefore, multiplying both sides of (\ref{factor}) by $q$, we
find that
$$\alpha_{j_1}\cdots\alpha_{j_R}e=h_1e_1+\dots+h_ke_k$$
with $h_i\in S(\Gs)$. Thus, in (\ref{factor}), none of
the weights $\alpha:G\rightarrow\Zt$ that divide the denominators
$p_i$ vanish on $\H$.  Hence, if $e\in \im(i_{\H}^*)$ for all $\H$,
then $p_i=1$ for all $i$, and so (\ref{factor}) tells us that
$$e=g_1e_1+\dots+g_ke_k\in\HGM$$
and the proof is complete.
\end{proofof}

Now suppose that $Z_{\H}$ is a connected component of $M^{\H}$ for some
subgroup $\H$ of $\G$ of order $|\H|=2^{n-1}$.  Let $i_{Z_{\H}}$ be the
inclusion
$$
i_{Z_{\H}}:Z_{\H}^{\G}\to Z_{\H}
$$
of the fixed points of $Z_{\H}$ into $Z_{\H}$.  Let $r_{Z_{\H}}$ be the
inclusion
$$
r_{Z_{\H}}:Z_{\H}^{\G}\to M^{\G}
$$
of the fixed points of $Z_{\H}$ into all of the fixed points.  Then we
have the following corollary of Theorem~\ref{Tchang}.

\begin{cor}\label{co:genCS}
Suppose that $\HGM$ is a free $\HG$-module.  A class
$$
f\in\HG(M^{\G};\Zt)
$$
is in the image of $i^*$ if and only if
$$
r_{Z_{\H}}^*(f)\in i_{Z_{\H}}^*(\HG(Z_{\H};\Zt))
$$
for every subgroup $\H$ of $\G$ of order $|\H|=2^{n-1}$ and every
connected component $Z_{\H}$ of $M^{\H}$.
\end{cor}

\begin{PROOF}
The proof is analogous to the proof of Theorem~$1$ in \cite{gh}.
It follows directly from Theorem~\ref{Tchang}.
\end{PROOF}

\section{Real GKM}\label{se:realGKM}

The goal of this section is to prove an analogue of
Theorem~\ref{th:GKM} for the real locus $X$ of $M$.  The proof will
require two hypotheses on $X$, namely
\begin{equation}\label{eq:fixed}
\XG=M^G
\end{equation}
and
\begin{equation}\label{eq:1skel}
X^{(1)}=X\cap M^{(1)},
\end{equation}
where $M^{(1)}$ is the one-skeleton of $M$ and $X^{(1)}$ the
one-skeleton of $X$.  We will begin by analyzing these conditions and
their implications.  We first note that the analogues for $M$ of the
conditions (\ref{eq:fixed}) and (\ref{eq:1skel}), namely
\begin{equation}
\# M^G<\infty
\end{equation}
and
\begin{equation}
\dim(M^{(1)})\leq 2,
\end{equation}
have a very simple and elegant interpretation in terms of the isotropy
representations of $G$ at fixed points of $M$.
\begin{thm}
\label{th:GKMhyp}
The conditions $\# M^G<\infty$ and $\dim(M^{(1)})\leq 2$ are satisfied
if and only if, for every $p\in M^G$, the weights $\tilde{\alpha}_{i,p}$,
$i=1,\dots, d$ of the isotropy representation of $G$ on $T_pM$ are
pair-wise linearly independent, that is for $i\neq j$, $\tilde{\alpha}_{i,p}$
is not a multiple of $\tilde{\alpha}_{j,p}$.
\end{thm}
For the proof of this, see \cite{GZ}.  When $M$ satisfies the two
conditions \eqref{eq:fixed} and \eqref{eq:1skel}, we say that $M$
is a {\em GKM manifold}. Let $\Z^*_G$ be the weight lattice of
$G$. By the $\!\!\mod 2$ reduction of a weight
$\tilde{\alpha}\in\Z^*_G$, we mean its image $\alpha$ in
$\Z^*_G/2\Z^*_G$. We will prove a real analogue of
Theorem~\ref{th:GKMhyp}.
\begin{thm}\label{th:realRules}
Suppose $M$ satisfies the hypotheses of Theorem~\ref{th:GKMhyp}.
Then the conditions $\XG=M^G$ and $X^{(1)}=X\cap M^{(1)}$ are
satisfied if and only if, for every $p\in M^G$, the $\!\!\mod 2$ reduced
weights, $\alpha_{i,p}$, are all distinct and non-zero.
\end{thm}
\begin{PROOF}
Let $Y$ be a connected component of $M^{\G}$.  Then $Y$ is a
$G$-invariant symplectic submanifold of $M$, and the action of $G$ on
it is Hamiltonian, so it contains at least one $G$-fixed point $p$.
However, the hypotheses above imply that the linear isotropy action of
$\G$ on $T_pM$ has no fixed points other than the origin.  Hence,
$\dim(Y)=0$ and $Y=\{ p\}$.  This argument applies to all the
connected components of $M^{\G}$, hence the connected components are
just the fixed points of $G$, and thus $\XG=M^G$.

The proof that $X^{(1)}=X\cap M^{(1)}$ is similar.  Let $\H$ be a
subgroup of $\G$ of index $2$, and let $Y$ be a connected
component of $M^{\H}$.  Then $Y$ is a $G$-invariant submanifold of
$M$, and because $\sigma\circ\tau_g  =  \tau_{g^{-1}}\circ\sigma$,
it is also $\sigma$-invariant.  Let $p\in Y$ be a $G$-fixed point, and
let
$$
T_pM=V_1\oplus\cdots\oplus V_d
$$
be the decomposition of $T_pM$ into the $2$-dimensional weight spaces
corresponding to the $\tilde{\alpha}_{i,p}$.  By the hypotheses on the
reduced weights $\alpha_{i,p}$, either
$$
(T_pM)^{\H}=\{ 0\},
$$
in which case $Y=\{ p\}$, or
\begin{equation}\label{eq:tang}
(T_pM)^{\H}=V_i=T_pY
\end{equation}
for some $i$.  Let $\chi_i$ be the character of $G$ associated with
the representation of $G$ on $V_i$ and let $H=\ker(\chi_i)$.  Then
$\H\subset H$ and
$$
(T_pM)^H=V_i.
$$
Thus, by (\ref{eq:tang}), $Y$ is the connected component of $M^H$
containing $p$, and in particular, $Y$ is contained in $M^{(1)}$.
Thus,
$$
Y^{\sigma}\subseteq X\cap M^{(1)}.
$$
Applying this argument to all index $2$ subgroups $\H$ of
$\G$ and all connected components of the fixed point sets of these
groups, one obtains the inclusion
$$
X^{(1)}\subseteq X\cap M^{(1)}.
$$
The reverse inclusion is obvious.  This completes the proof.
\end{PROOF}

The hypotheses of Theorem~\ref{th:realRules} impose some rather severe
restrictions on the manifold $M$.  For instance, the cardinality of
the set of $\!\!\mod 2$ reduced weights, $\Z^*_G/2\Z^*_G$, is $2^n$.
Therefore, since the reduced weights $\alpha_{i,p}$ are distinct
and non-zero for $i=1,\dots,d$, we must have that $d\leq 2^n-1$.
Hence,
\begin{equation}\label{eq:dim}
\dim(M)=2d\leq 2^{n+1}-2.
\end{equation}
For example, if $n=2$, then $\dim(M)\leq 6$.  This leads us to make
the following definition.

\begin{definition}
If $M$ is a GKM manifold, and if for every $p\in
M^G$, the $\!\!\mod 2$ reduced weights, $\alpha_{i,p}$, are all
distinct and non-zero, we will say that $M$ is a {\em $\!\!\mod 2$ GKM
manifold}.
\end{definition}

Next, we show that relatively few compact {\em homogeneous}
symplectic manifolds (e.g. coadjoint orbits) are $\!\!\mod 2$ GKM
manifolds.  Consider coadjoint orbits of the classical compact
simple Lie groups associated with the Dynkin diagrams $A_n$,
$B_n$, $C_n$ and $D_n$.  Let $\varepsilon_i$, for $i=1,\dots,n$,
be the standard basis vectors of $\R^n$.  The positive roots
associated to the Dynkin diagram $A_n$ consist of
$$
\varepsilon_i-\varepsilon_j,\ \  i<j;
$$
so their $\!\!\mod 2$ reductions are distinct and non-zero.
However, for $B_n$, $C_n$, and $D_n$, this list of positive roots
contains
$$
\varepsilon_i-\varepsilon_j\mbox{ and } \varepsilon_i+\varepsilon_j,\
\  i<j,
$$
so we conclude
\begin{thm}
Each coadjoint orbit of $SU(n)$ is a $\!\!\mod 2$ GKM space.  However, for
other compact simple Lie groups, {\em no} maximal coadjoint orbit can
be a $\!\!\mod 2$ GKM space.
\end{thm}
On the other hand, on a more positive note, one has
\begin{thm}
If $M$ is a non-singular projective toric variety, then it is a
$\!\!\mod 2$ GKM space.
\end{thm}
\begin{PROOF}
If $M$ is a non-singular toric variety, the weights $\tilde{\alpha}_{i,p}$,
$i=1,\dots,n$, are a $\Z$-basis for $\Z^*_G$, so their images in
$\Z^*_G/2\Z^*_G$ are a $\Zt$ basis of $\Z^*_G/2\Z^*_G$.
\end{PROOF}
This theorem, combined with Theorem~\ref{th:realGKM}, gives us a
new description of the equivariant cohomology of a real toric
variety.  The ordinary and $\G$-equivariant cohomology of these real
loci has been computed by Davis and Januszkiewicz \cite{DJ}.  Their
description of these rings is analogous to Danilov's description of
the ordinary and $G$-equivariant cohomology of the original toric
varieties.

We will now prove a real locus version of the GKM theorem with $\Zt$
coefficients.  Recall from Section~\ref{se:intro} that Theorem~\ref{th:GKM}
of GKM characterizes the image of $i^*:H_G^*(M;\C)\to H_G^*(M^G;\C)$
in terms of the weights of the isotropy representations of $G$ on the tangent
spaces at the fixed points.

To prove an analogue of this for the real locus of a symplectic
manifold, we must first compute the $\Zt$-equivariant cohomology with
$\Zt$ coefficients of $\R P^1$.  Recall that $S^1$ acts on $\C P^1$ by
$\theta\cdot [z_0:z_1]=[z_o:e^{i\theta}z_1]$.  This is a Hamiltonian
action, with respect to the Fubini-Study symplectic form on $\C P^1$.
Furthermore, complex conjugation is an anti-symplectic involution on
$\C P^1$, with fixed point set $\R P^1$.  There is a residual action
of $\Zt$ on $\R P^1\iso S^1$ which reflects $S^1$ about the $y$-axis.

\begin{lemma}\label{rp1}
Let $N$ and $S$ denote the fixed points of the $\Zt$ action on $\R
P^1$. Then the image of the map
$$
i^*:H_{\Zt}^*(\R P^1;\Zt)\to H_{\Zt}^*(N;\Zt)\oplus H_{\Zt}^*(S;\Zt)
$$
is the set of pairs $(f_N, f_S)$ such that
$$
f_N+f_S\in x\cdot\Zt[x].
$$
\end{lemma}

\begin{PROOF}
It is clear that the constant functions are equivariant classes in
$$
H_{\Zt}^0(\R P^1;\Zt).
$$
Furthermore, we know that $\dim H_{\Zt}^0(\R P^1;\Zt)=1$, and so these
are the only equivariant classes.  Finally, $\dim H_{\Zt}^i(\R
P^1;\Zt)=2$ for $i>0$, and so indeed, the condition stated is the only
condition of pairs $(f_N,f_S)\in H_{\Zt}^*(N;\Zt)\oplus
H_{\Zt}^*(S;\Zt)$.
\end{PROOF}

Theorem~\ref{th:realGKM} identifies the image of the map
$$
i^*:\HG(X;\Zt)\to\HG(X^{\G};\Zt)
$$
in terms of weights of isotropy representations of $\G$ on the tangent
spaces at the fixed points.

\begin{proofof}{Proof of Theorem~\ref{th:realGKM}}
The result follows immediately from Corollary~\ref{co:genCS} and
Lemma~\ref{rp1}.
\end{proofof}

The results of this section and the previous section have been proved
independently by Schmid \cite{schmid}.  Schmid uses an equivariant
Morse theoretic approach, and consequently the proofs are quite different.

As a result of equivariant formality, we get two corollaries of
Theorem~\ref{th:realGKM} concerning the relation between the ring
structure of the cohomology of $M$ and the cohomology of $X$.

\begin{cor}\label{co:ring1}
Suppose that $M$ is a GKM manifold and a $\!\!\mod 2$ GKM manifold.  Then
there is a graded ring isomorphism
$$
H_G^{2*}(M;\Zt)\iso H^*_{\G}(X;\Zt).
$$
\end{cor}

\begin{cor}\label{co:ring2}
Suppose that $M$ is a GKM manifold and a $\!\!\mod 2$ GKM manifold.  Then
there is a graded ring isomorphism
$$
H^{2*}(M;\Zt)\iso H^*(X;\Zt).
$$
\end{cor}

Note that this last corollary strengthens Duistermaat's original
result from an isomorphism of vector spaces to an isomorphism of
rings.

\begin{remark}
Some of the results of this section, most importantly
Theorem~\ref{th:realGKM}, are valid not only for the real locus
$X$ of a Hamiltonian $G$-manifold, but more generally for any
compact $\G$-manifold $X$ which satisfies the following
properties:
\begin{enumerate}
\item $X$ is equivariantly formal;
\item $X^{\G}$ is finite; and
\item the weights of $X$ satisfy the properties of a $\!\!\!\mod 2$ GKM
manifold.
\end{enumerate}
In this situation, we may still characterize the structure of the
one-skeleton. Theorem~\ref{th:realGKM} still follows from
injectivity and the Chang-Skjelbred theorem.
\end{remark}

\section{Real GH}\label{se:realGH}

Goldin and Holm generalize Theorem~\ref{th:GKM} to the case where the
one-skeleton has dimension at most $4$.  The goal of this section is
to prove a real version of the Goldin-Holm theorem with $\Zt$
coefficients.  Again, we require the hypotheses that the
$(\Zt)^n$-fixed points of the real locus are the same as the $G$-fixed 
points of $M$ as in (\ref{eq:fixed}); and that the real locus of the
one-skeleton is the same as the one-skeleton of the real locus, as in
(\ref{eq:1skel}). Finally, we require
$$
\# M^G<\infty,
$$
and
$$
\dim(M^{(1)})\leq 4.
$$
If a manifold satisfies these last two hypotheses, we will say that it
is a {\em GH manifold}. These hypotheses have a nice interpretation in
terms of the isotropy representations of $G$ at the fixed points of
$M$.

\begin{theorem}\label{th:GHrealRules}
The conditions $\# M^G<\infty$ and $\dim(M^{(1)})\leq 4$ are satisfied 
if and only if the weights $\alpha_{i,p}$ of the isotropy
representation of $G$ on $T_pM$ have the property that every three
span a vector subspace of dimension at least two.
\end{theorem}

These hypotheses on $M$ have real analogues, namely that $\# X<\infty$
and the one-skeleton $X^{(1)}$ of $X$ is at most $2$-dimensional.
We will state without proof the following real analogue of
Theorem~\ref{th:GHrealRules}.

\begin{theorem}
Suppose that $M$ satisfies the hypotheses of
Theorem~\ref{th:GHrealRules}. If the conditions $X^{\G}=M^G$ and
$X^{(1)}=M^{(1)}\cap X$ are satisfied, then for every $p\in M^G$, the
$\!\!\mod 2$ reduced weights $\alpha_{i,p}^{\#}$ are all non-zero, and
each element of $S((\G)^*)=\Zt[x_1,\dots,x_n]$ appears no more than
twice.
\end{theorem}

The proof of this theorem is nearly identical to that of
Theorem~\ref{th:realRules}.  The hypotheses of this theorem, although
weaker than those of Theorem~\ref{th:realRules}, still impose
restrictions on the manifold $M$.  The cardinality of the set of $\!\!\mod 
2$ reduced weights is $2^n$.  Since the weights are non-zero, and each 
weight can appear at most twice,
$$
d\leq 2\cdot (2^n-1),
$$
and so
$$
\dim (M)=2d\leq 2\cdot(2\cdot(2^n-1))=2^{n+2}-4.
$$
For instance, if $n=2$, $\dim(M)\leq 12$.  We will now show an example 
where the condition that the reduced weights be non-zero is {\em not}
satisfied.

\begin{Eg}{
Consider $\C P^2$ with homogeneous coordinates
$[z_0:z_1:z_2]$.  Let $T=S^1$ act on $\C P^2$ by 
$$
e^{i\theta}\cdot [z_0:z_1:z_2] = [e^{-i\theta}z_0:z_1:e^{i\theta}z_2].
$$
This action has three fixed points: $[1:0:0]$, $[0:1:0]$ and
$[0:0:1]$.  

The weights at these fixed points are
$$
\begin{array}{cc}
\mbox{Fixed point} & \mbox{Weights} \\
\mbox{$p_1=[1:0:0]$} & x, 2x \\
\mbox{$p_2=[0:1:0]$} & -x, x \\
\mbox{$p_3=[0:0:1]$} & -2x, -x
\end{array}
$$ 
where we have identified $\mathfrak{t}^*$ with degree one
polynomials in $\C [x]$. As cohomology elements, these are assigned
degree two. Using Theorem~\ref{th:4total} below, we can
compute the $S^1$ equivariant cohomology of $\C P^2$ as follows.
The image of the equivariant cohomology $H_{S^1}^*(\C P^2
)$ in 
$$
H_{S^1}^*(\{ p_1,p_2,p_3\} )\iso \bigoplus_{i=1}^3\C [x]
$$
is the subalgebra generated by the triples of functions
$(f_1,f_2,f_3)$ such that 
$$
f_i-f_j\in x\cdot\C [x]  \mbox{ for every $i$ and $j$, and } 
$$
$$
\frac{f_1}{2x^2}-\frac{f_2}{x^2}+\frac{f_3}{2x^2}\in\C [x]. 
$$
However, when we try to compute the $\Zt$ equivariant cohomology of
$\R P^2$, the real locus of $\C P^2$, we run into a problem.  The
$\!\!\mod 2$ reduced weights are given in the table below.
$$
\begin{array}{cc}
\mbox{Fixed point} & \mbox{Weights} \\
\mbox{$p_1=[1:0:0]$} & x, 0, \\
\mbox{$p_2=[0:1:0]$} & x, x, \\
\mbox{$p_3=[0:0:1]$} & 0, x.
\end{array}
$$
The problem with this $\Zt$ action on $\R P^2$ is that it no longer
has isolated fixed points.  There is an entire $\R P^1$ which is fixed 
by this $\Zt$ action.
Thus, we cannot compute the $\Zt$ equivariant cohomology of $\R P^2$
using these methods. 
}\end{Eg}

We make the following definition, analogous to the definition of $\!\!\mod
2$ GKM manifolds given in Section~\ref{se:realGKM}.

\begin{definition}
Suppose that $M$ is a GH manifold, and furthermore that $X^{\G}=M^G$ and
$X^{(1)}=M^{(1)}\cap X$.  In this case, we will say that $M$ is a {\em
$\!\!\mod 2$ GH space}. 
\end{definition}

Recall the following
properties about the $G$-equivariant cohomology of manifolds with
one-skeleta of dimension at most $4$.  These are proved in \cite{gh},
although the reader is cautioned to the different notation used there.
First, we compute the $S^1$ equivariant cohomology of a $4$-manifold,
and then we use this computation to determine the equivariant
cohomology of any manifold with one skeleton of dimension at most $4$.

\begin{lemma}
\label{pr:GKMin4} Let $X$ be a compact, connected symplectic
4-manifold with an effective Hamiltonian $S^1$ action with isolated
fixed points $X^{S^1}=\{ p_1,\dots ,p_d\}$.
The map $i^*:H_{S^1}^*(X)\rightarrow H_{S^1}^*(X^{S^1})$ induced by
inclusion is an injection with image
$$
\left\{ (f_1,\dots ,f_d)\in\bigoplus_{i=1}^d S(\algs^*)\ \left| \ f_i-f_j\in
x\cdot\C [x],\ \sum_{i=1}^d\frac{f_i}{\alpha_1^i\alpha_2^i}\in
S(\algs^*)\right.\right\},
$$
where $\alpha_1^i$ and $\alpha_2^i$ are the (linearly
dependent) weights of the $S=S^1$ isotropy action on $T_{p_i}X$.
\end{lemma}

\begin{thm}
\label{th:4total}
Let $M$ be a
compact, connected symplectic manifold with an effective Hamiltonian
$G$-action.  Suppose further that the $G$-action has only isolated
fixed points $M^G=\{ p_1,\dots ,p_d\}$ and that the one skeleton has
dimension at most $4$. Let $f_i\in H_G^*$ denote the restriction
of $f\in H_G^*(M)$ to the fixed point $p_i$.
The image of the injection $i^*:H_{G}^*(M)\rightarrow H_{G}^*(M^G)$ is 
the subalgebra of functions $(f_1,\dots
,f_d)\in\bigoplus_{i=1}^d  S(\g^*)$ which satisfy
$$
\left\{\begin{array}{ll}
\pi_H^*(f_{i_j})=\pi_H^*(f_{i_k}) & \mbox {if } \{ p_{i_1},\dots
,p_{i_{\ell}}\}=Z_H^G  \\
\sum_{j=1}^l\frac{f_{i_j}}{\alpha_1^{i_j}\alpha_2^{i_j}}\in
S(\g^*) & \mbox {if }\{ p_{i_1},\dots ,p_{i_{\ell}}\}=Z_H^G \mbox{ and }
\dim Z_H =4
\end{array}
\right . 
$$
for all $H\subset G$ codimension-1 tori and all connected components
$Z_H$ of $M^H$, where $\alpha_1^{i_j}$;
$\alpha_2^{i_j}$ are the (linearly dependent) weights of the $G$
action on $T_{p_{i_j}}Z_H$; and $\pi_H:\algh\into \algg$ is inclusion.
\end{thm}

We can use these computations to compute the $(\Zt)^n$ equivariant
cohomology of a $\!\!\mod 2$ GH manifold.

\begin{lemma}
\label{le:realGKMin4} Let $M$ be a compact, connected symplectic
4-manifold with an effective Hamiltonian $S^1$ action with isolated
fixed points $M^{S^1}=\{ p_1,\dots ,p_d\}$.  Suppose further that $M$
is a $\!\!\mod 2$ GH manifold with real locus $X$. The map
$$
i^*:H_{\Zt}^*(X;\Zt)\rightarrow 
H_{\Zt}^*(X^{\Zt};\Zt)
$$
induced by inclusion is an injection with image
\begin{equation}
\label{eq:ZtonM^4}
\left\{ (f_1,\dots ,f_d)\in\bigoplus_{i=1}^d \Zt[x]\ \left|
\begin{array}{l}\ f_i-f_j\in x\cdot\Zt [x],\\
\sum_{i=1}^d\frac{f_i}{\alpha_1^i\alpha_2^i}\in 
\Zt[x]\end{array}\right.\right\},
\end{equation}
where $\alpha_1^i$ and $\alpha_2^i$ are the linearly
dependent weights of the $\Zt$ isotropy representation on $T_{p_i}X$. (In
this case, $\alpha_1^i=\alpha_2^i=x$.)
\end{lemma}

\begin{PROOF}
The map $i^*$ is injective because $X$ is equivariantly formal.  We
know that the $f_i$ must satisfy the first condition because the
functions constant on all the vertices are the only equivariant
classes in degree $0$, as $\dim H_{\Zt}^0(X;\Zt)=1$.  The second
condition is necessary as a direct result of the $\Zt$ version of the
localization theorem proved in Section~\ref{se:cs}.   Notice that this
condition gives us one relation in degree $1$ cohomology.  A dimension
count shows us that these conditions are sufficient.  As an
$S((\Zt)^*)$-module, $H_{\Zt}^*(X;\Zt)\iso 
H^*(X;\Zt)\otimes H_{\Zt}^*(pt;\Zt)$.  Thus, the equivariant
Poincar\'{e} polynomial is
\begin{eqnarray*}
P^{\Zt}_t(X) & = & (1+(d-2)t+t^2)\cdot (1+t+t^2+\dots) \\
         & = & 1+(d-1)t+dt^2+\dots+dt^{n}+\cdots.
\end{eqnarray*}
\noindent  As $H_{\Zt}^*(X;\Zt)$ is generated in degree 1, the $d-1$
degree 1 classes given by the $(f_1,\dots,f_d)$ subject to the localization
condition generate the
entire cohomology ring. Thus we have found all the conditions.
\end{PROOF}

We now prove Theorem~\ref{th:realGH}, computing the cohomology of any
$\!\!\mod 2$ GH manifold.  We will show that the image of
$i^*:H_{\G}^*(X)\rightarrow H_{\G}^*(X^{\G})$ is the subalgebra of 
functions $(f_1,\dots ,f_d)\in\bigoplus_{i=1}^d S(\G^*)$ which satisfy
$$
\left\{\begin{array}{ll}
\pi_{\H}^*(f_{i_j})=\pi_{\H}^*(f_{i_k}) & \mbox {if } \{ p_{i_1},\dots
,p_{i_l}\}=Z_{\H}^{\G}  \\
\sum_{j=1}^l\frac{f_{i_j}}{\alpha_1^{i_j}\alpha_2^{i_j}}\in
S(\G^*) & \mbox {if }\{ p_{i_1},\dots ,p_{i_l}\}=Z_{\H}^G \mbox{ and }
\dim Z_{\H} =4,
\end{array}
\right . 
$$
for all subgroups $\H$ of $\G$ of order $|\H|=2^{n-1}$ and all
connected components $Z_{\H}$ of $X^{\H}$.

\begin{proofof}{Proof of Theorem~\ref{th:realGH}}
This follows immediately from Corollary~\ref{co:genCS} and
Lemma~\ref{le:realGKMin4}.
\end{proofof}

There are two immediate corollaries in this setting, analogous to
Corollaries~\ref{co:ring1} and \ref{co:ring2}.

\begin{cor}\label{co:ring3}
Suppose that $M$ is a GH manifold,  and that $M^G=X^{\G}$ and
$M^{(1)}\cap X=X^{(1)}$.   Then there is a graded ring isomorphism
$$
H_G^{2*}(M;\Zt)\iso H^*_{\G}(X;\Zt).
$$
\end{cor}

\begin{cor}\label{co:ring4}
Suppose that $M$ is a GH manifold,  and that $M^G=X^{\G}$ and
$M^{(1)}\cap X=X^{(1)}$.  Then there is a graded ring isomorphism
$$
H^{2*}(M;\Zt)\iso H^*(X;\Zt).
$$
\end{cor}

\section{Applications to String Theory}\label{se:strings}

Consider the $\Zt$ action on $T^n$, which reflects each copy of $S^1$.  
Then the equivariant cohomology ring
$$
H^*_{\Zt}(T^n;\Zt)
$$
classifies all possible orientifold configurations of Type II string
theories, compactified on $T^n$. See Section 3 and Appendix C of
\cite{string} for more details. Yang-Hui He pointed this example out
to us. Using the results of Section~\ref{se:realGKM}, we can now
compute this equivariant cohomology.

First, we recognize $T^n$ as the real locus of $M=\C
P^1\times\cdots\times\C P^1=(\C P^1)^n$.  This space $M$ has a natural 
$T^n$ action, where the $i$th copy of $S^1$ acts in the
standard fashion on the $i$th copy of $\C P^1$.  We can
compute the $(\Zt)^n$-equivariant cohomology of this space quite
easily.  The GKM graph associated to $(\C P^1)^n$ with the $T^n$
action described above is the $n$-dimensional hypercube.  The
vertices correspond to the binary words of length $n$.  
Two binary words are connected by an edge if they differ in exactly
one bit.  Suppose $v$ and $w$ differ in exactly the $i$th
bit.  Then the weight associated to the edge $(v,w)$ is $x_i$.  Thus,
when $n=3$, the GKM graph and weights are shown in the figure below.

\begin{figure}[h]
\centerline{
\epsfig{figure=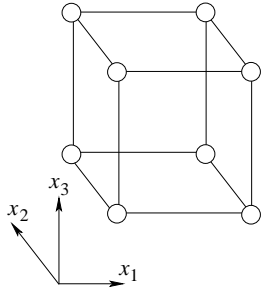,height=1.5in}
}
\smallskip
\centerline{
\parbox{4.5in}{\caption[hypercube]{{\small This shows the GKM graph
and the weights for $(\C P^1)^3$.}}}
}
\end{figure}

Notice that the reduced weights are all non-zero and are distinct in
$\Z_G/2\Z_G$.  Thus, we can apply Theorem~\ref{th:realGKM} to compute
$$
H_{(\Zt)^n}(T^n;\Zt).
$$
That is, the equivariant cohomology is the set of functions
$f:V\to\Zt[x_1,\dots,x_n]$ such that for every edge $(v,w)\in E$, we
have
$$
f(v)+f(w)\in x_i\cdot\Zt[x_1,\dots,x_n].
$$

We can now consider the copy of $\Zt$ sitting diagonally inside
$(\Zt)^n$.  This copy of $\Zt$ acts on $T^n$, and this is the action
that originally interested physicists.  We can now compute 
the $\Zt$-equivariant cohomology simply by projecting
$$
\pi:S(((\Zt)^n)^*)=\Zt[x_1,\dots,x_n]\to \Zt[x]=S((\Zt)^*)
$$
where $x_i$ gets sent to $x$.  Then
$$
H_{\Zt}^*(T^n;\Zt)=\pi(H_{(\Zt)^n}(T^n;\Zt)).
$$

\section{Recreational applications}\label{se:recreation}

We will describe below a (somewhat idealized) real world application
of the results of Section~\ref{se:realGKM}. Let $S=\{1,\dots,n\}$ be a 
list of companies whose stocks are being traded in the stock market.
Let $A_i$ for $i=1,\dots,d$ be investors, and let $S_i\subseteq S$ be
the portfolio of $A_i$.  Suppose that for certain pairs of investors
$A_i$ and $A_j$, the symmetric difference
\begin{equation}\label{eq:8.1}
(S_i-S_j)\cup(S_j-S_i)=S_e,\ \ e=(i,j),
\end{equation}
i.e. the relative status of the portfolios of $A_i$ and $A_j$, is
given.  To what extent does this information determine the portfolios
$S_i$? It is easy to see that it cannot {\em uniquely} determine the
$S_i$'s.  If $S_1,\dots,S_n$ is one solution to (\ref{eq:8.1}), one gets 
another solution by taking a fixed subset $S_0$ of $S$ and replacing
the $S_i$'s by the symmetric difference
\begin{equation}\label{eq:8.2}
S_i'=(S_i-S_0)\cup(S_0-S_i).
\end{equation}
Therefore, we will slightly rephrase this question.  Let $E$ be a
collection of two-element subsets of $S$, and suppose that for every
$e\in E$, one is given a subset $S_e$ of $S$.  List all solutions
$S_1,\dots,S_n$ of the Boolean identities (\ref{eq:8.1}).  We will say
two solutions are identical if they satisfy (\ref{eq:8.2}).  Note, by
the way, that (\ref{eq:8.2}) can be rewritten
$$
S_i=((S_i'-S_0)\cup S_0)-S_i'.
$$
Hence, we have defined an equivalence relation.

One can inject an element of randomness into this problem by positing
that, for $e=(i,j)\in E$,
\begin{equation}\label{eq:8.3}
(S_i-S_j)\cup(S_j-S_i)  \in  \{	S_e ,\varnothing\}.
\end{equation}
In other words, either the symmetric difference is given by
(\ref{eq:8.1}) or $S_i=S_j$. Again, the problem is to list all
possibilities for the $S_i$'s.  Clearly, the solutions of (\ref{eq:8.3})
contain the solutions to (\ref{eq:8.1}); so by solving (\ref{eq:8.3}), one
gets an upper bound on the number of solutions to (\ref{eq:8.1}).
Moreover, there are a lot of trivial solutions of (\ref{eq:8.3}),
namely,
\begin{equation}\label{eq:8.4}
S_i=S_0,\ i=1,\dots,n,
\end{equation}
where $S_0$ is, as above, a fixed subset of $S$. These can immediately
be discarded as potential solutions of (\ref{eq:8.1}).

There is an elegant way of reformulating (\ref{eq:8.1}) and (\ref{eq:8.3})
in the language of $\mod 2$ arithmetic.  Let $\Gamma$ be the graph
whose vertices are the $A_i$'s and whose edges are the members of
$E$.  For every edge $e\in E$, let $\alpha_e$ be the element of
$\Zt^n$ whose $k$th coordinate is $1$ if and only if $k\in
S_e$ and $\alpha_i$ the element of $\Zt^n$ whose $k$th
coordinate is $1$ if and only if $k\in S_i$.  Then (\ref{eq:8.1}) is
equivalent to
\begin{equation}\label{eq:8.5}
\alpha_i+\alpha_j=\alpha_e
\end{equation}
and (\ref{eq:8.3}) is equivalent to
\begin{equation}\label{eq:8.6}
\alpha_i+\alpha_j=\lambda\alpha_e,\ \lambda\in\Zt.
\end{equation}
In particular, (\ref{eq:8.3}) becomes an identity of the type described
in Theorem~\ref{th:realGKM}, the multiple of $\alpha_e$ on the right
being in the degree zero component of $\Zt[x_1,\dots,x_n]$.

Now let $M$ be a GKM manifold acted on by an $n$-torus $T$, and let
$\Gamma$ be its associated graph.  If the $\alpha_e$'s in
(\ref{eq:8.6}) are the weights of $\G$ which were defined in
Theorem~\ref{th:realGKM} the results of Section~\ref{se:realGKM} tell
us that solutions of (\ref{eq:8.6}) can be identified with elements of
$$
H^1_{\G}(M^{\sigma};\Zt).
$$
We will use this observation to determine the solutions of
(\ref{eq:8.6}) in a couple simple, but interesting, examples.

\begin{Eg}{
{\bf The complete graph on $n$ vertices.}  Consider $M=\C P^{n-1}$ as a
$T^n$ manifold (ignoring the fact that the diagonal subgroup of $T^n$
acts trivially).  The corresponding graph is the complete graph on $n$
vertices: its vertices are $A_1,\dots,A_n$ and every pair of vertices
is joined by an edge.  The weights $\alpha_{(i,j)}$ from
Theorem~\ref{th:realGKM} are just $x_i+x_j$, so the conditions
(\ref{eq:8.1}) become
\begin{equation}\label{eq:8.7}
(S_i-S_j)\cup(S_j-S_i)=\{ i,j\}.
\end{equation}
The solutions of (\ref{eq:8.7}) are in one to one correspondence with
the elements of
$$
H^1_{\G}(\R P^{n-1};\Zt)
$$
or, alternatively, of
\begin{equation}\label{eq:8.8}
(H^0(\R P^{n-1};\Zt)\otimes\Zt^n)\oplus H^1(\R P^{n-1};\Zt).
\end{equation}
The elements of the first summand correspond to the trivial solutions
of (\ref{eq:8.7}).  So if we identify solutions which are equivalent in
the sense of (\ref{eq:8.2}), the non-trivial solutions of (\ref{eq:8.7})
correspond to non-zero elements of $H^1(\R P^{n-1};\Zt)$.  However,
$H^1(\R P^{n-1};\Zt)=\Zt$, so there is just one non-trivial solution
up to equivalence, and it is given by
$$
S_i=\{ i\};
$$
this is also the unique solution of (\ref{eq:8.1}) up to equivalence.
}
\end{Eg}

\begin{Eg}{
{\bf The permutahedron.}
Let $M$ be the complex flag variety $U(n)/T^n$, and consider the $T^n$
action on $M$ by left multiplication.  The graph associated with $M$
is the {\em permutahedron}.  Its vertices are the elements of the
symmetric group, $S_n$, and two vertices $\sigma$ and $\tau$ are
joined by and edge if $\tau\sigma^{-1}$ is a transposition.  If $e$ is
the edge joining $\sigma$ to $\tau$ and $\tau\sigma^{-1}$ is the
transposition switching $i$ and $j$, then as in the previous example,
$\alpha_e=x_i+x_j$, so the conditions (\ref{eq:8.1}) become
$$
(S_{\sigma}-S_{\tau})\cup(S_{\tau}-S_{\sigma})=\{ i,j\}.
$$
As in the previous example, the non-trivial solutions to (\ref{eq:8.3})
can be identified with the non-zero elements of $H^1(M^{\sigma};\Zt)$,
and since $M^{\sigma}$ is the real flag variety,
$$
H^1(M^{\sigma};\Zt)\iso \Zt^{n-1}
$$
If one thinks of $\Zt^{n-1}$ as the quotient of $\Zt^n$ by the
diagonal subgroup $(\lambda,\dots,\lambda),$ with $\lambda\in \Zt$,
the solutions corresponding to $\alpha\in\Zt^n\mod
(\lambda,\dots,\lambda)$ is given by 
\begin{equation}\label{eq:8.11}
S_{\sigma}=\sigma(S_0),
\end{equation}
where $k\in S_0$ if and only if the $k$th coordinate of
$\alpha$ is $1$.  Notice, by they way, if we replace $\alpha$ by
$\alpha +(1,\dots,1)$, (\ref{eq:8.11}) becomes
\begin{equation}\label{eq:8.12}
S_{\sigma}^c=\sigma(S_0^c),
\end{equation}
where $S_{\sigma}^c$ and $S_0^c$ are the complements of $S_{\sigma}$
and $S_0$ in $S$, so (\ref{eq:8.11}) and (\ref{eq:8.12}) are equivalent.
It is easy to see by inspection that none of the solutions
(\ref{eq:8.11}) of (\ref{eq:8.3}) are also solutions of (\ref{eq:8.1}).

An interesting special case of this example is the complete bipartite
graph $K_{3,3}$.  In this case, the vertices of $\Gamma$ are $A_1$,
$A_2$, $A_3$ and $B_1$, $B_2$, and $B_3$; the edges are all pairs
$(A_i,B_j)$; and the sets $S_e$ are
$$
S_{A_i,B_j}=\left\{\begin{array}{ll}
			\{ i,j\} & i\neq j\\
			\{ 1,2,3\}-\{ i\} & i=j
		   \end{array}\right.
$$
In this example, (\ref{eq:8.3}) has one non-trivial solution, up to
equivalence, namely $S_{A_i}=S_{B_i}=\{ i\}$.
}\end{Eg}

\end{document}